\documentclass[10pt,reqno]{article}
\usepackage[letterpaper,hmargin=1in,top=1in,bottom=1.2in,footskip=0.6in]{geometry}
\usepackage[OT1]{fontenc}
\usepackage{amsthm,amsmath,amssymb,mathtools}
\usepackage[numbers]{natbib}
\usepackage{enumitem}
\usepackage{mathrsfs,bm}
\usepackage{xcolor}
\usepackage[colorlinks,linkcolor=blue,citecolor=blue,urlcolor=blue]{hyperref}
\hypersetup{pdftitle={LFL for Gaussian beta corner process}}

\numberwithin{equation}{section}
\theoremstyle{plain}
\newtheorem{thm}{Theorem}[section]

\newtheorem{cor}[thm]{Corollary}
\newtheorem{pro}[thm]{Proposition}

\newtheorem{defn}[thm]{Definition}
\theoremstyle{remark}

\theoremstyle{plain}
\newtheorem*{claim*}{Claim}

\newcommand{\E}{{\mathbb E}}

\renewcommand{\P}{{\mathbb P}}

\newcommand{\dd}{\mathrm{d}}
\newcommand{\bs}{\boldsymbol}
\renewcommand{\mathbf}[1]{\bs{#1}}

\renewcommand{\le}{\leqslant}
\renewcommand{\geq}{\geqslant}
\renewcommand{\ge}{\geqslant}
\renewcommand{\epsilon}{\varepsilon}

\newcommand{\cC}{\mathcal C}
\newcommand{\cF}{\mathcal F}
\newcommand{\one}{\mathbf 1}

\newcommand{\rR}{{\mathbf{c}}_{\mathrm R}}
\newcommand{\rL}{{\mathbf{c}}_{\mathrm L}}

\DeclareMathOperator{\GammaLaw}{Gamma}
\DeclareMathOperator{\Dir}{Dir}

\allowdisplaybreaks[2]
\setlist[enumerate]{label=\textup{(\roman*)},leftmargin=*,itemsep=3pt}

\makeatletter
\def\paragraph{\@startsection{paragraph}{4}
  \z@\z@{-\fontdimen2\font}
  {\normalfont\bfseries}}
\makeatother

\title{Secular Dynamics and the LFL for $\beta$ Corner Process}

\author{
Zhigang Bao\thanks{
Department of Mathematics,
The University of Hong Kong.
Email: \texttt{zgbao@hku.hk}.}
\and
Ruizhe Xu\thanks{
Department of Mathematics,
The University of Hong Kong.
Email: \texttt{ruizhexu.math@gmail.com}.}
}

\date{}

\begin{document}

\maketitle

\begin{abstract} In this paper, we consider the Gaussian $\beta$ corner process, whose edge particles are known to  converge to the Tracy-Widom $\beta$ law. We study the almost sure behaviour of the extrema of the whole trajectory of the edge particles and establish a law of fractional logarithm (LFL) for the corner process for every fixed $\beta>0$, extending the previous known results in the case $\beta=1,2$.  Our proof relies on the upward Gamma and downward Dirichlet
secular equations, which can be regarded as two discrete dynamics to form the corner processes. The secular dynamics are used in both the correlation estimate and the decorrelation estimate, which are known as two key steps of proving LFL in \cite{BCEHK-LFL} for the Wigner minor process. Especially, in the decorrelation part, they are used as discrete analogues of the continuous Dyson Brownian motion. 
\end{abstract}

\medskip
\noindent\textbf{2020 Mathematics Subject
Classification.}
60B20, 60F15, 60J05.

\smallskip
\noindent\textbf{Keywords.}
Gaussian $\beta$ corner process,
law of fractional logarithm,
secular equation dynamics.

\section{Introduction and main results}\label{sec:intro}
\subsection{Motivation}\label{subsec:motivation}
There has been a vast of literature on the fluctuation of the largest eigenvalue of random matrix models. The most fundamental result is arguably the Tracy-Widom law, for GOE and GUE \cite{TW1994,TW1996}. There are two main directions to extend the result. One direction is to start with the matrix  model, and consider general real or complex Wigner matrices, and the other is to start from the jpdf of the eigenvalues and consider more general Gaussian $\beta$ ensemble with $\beta=1$ for GOE and $\beta=2$ for GUE. For the former, depending on the symmetry class, now it has been well known that the Tracy-Widom laws are universally true for general Wigner matrices \cite{Soshnikov1999,TaoVu2010Edge,Johansson2012Edge,LeeYin2014, ErdosYauYin2012}. For the latter, it is known from \cite{RRV} that the largest eigenvalue converges to the Tracy-Widom $\beta$ law.  If we view the Tracy-Widom law as an analogue of the CLT for the classical iid sum, it is natural to ask if there is an analogue of the classical law of iterated logarithm (LIL) ~\cite{HW}, for the sequence of the largest eigenvalues, if we consider a random matrix array. Such a question was explicitly asked by Kalai \cite{KalaiMO2013}. More specifically, if one starts with an infinite dimensional GOE or GUE, we can consider a sequence of nested principal minors of it.  Denote by $x_1^{(n)}$ as the largest eigenvalue of the $n$ by $n$ dimensional minor, and define the rescaled particles
\begin{equation}\label{eq:edgevariable}
  X_n:=n^{1/6}\bigl(x_1^{(n)}-2\sqrt n\bigr)
  \xrightarrow{d}\mathrm{TW}_\beta,
\end{equation}
for which we can ask what the almost sure extreme fluctuation scale of the whole sequence of $\{X_n\}_{n \geq 1}$ is. Similarly to the classical LIL, such problem does not only require us to understand  one individual $X_n$, it also requires one to investigate the correlation structure of the whole sequence. 

A pioneering work was done by  Paquette and
Zeitouni~\cite{PZ}, who initiated the study of this problem 
for the GUE minor process. They established
the sharp upper limit for $X_n$ and conjectured the corresponding
lower limit, which was only proved recently by
Baslingker, Basu, Bhattacharjee and
Krishnapur~\cite{BBBK-LFL}, again for GUE. More recently, Bao, Cipolloni,
Erd\H{o}s, Henheik and Kolupaiev~\cite{BCEHK-LFL}
proved a full version of LFL for all real symmetric and complex
Hermitian Wigner minor processes, using a robust martingale and dynamic approach which do not rely on any specific structure of GUE mnior process. Together, these works
establish, for the classical symmetry parameters
$\beta=1,2$, the following LFL
\begin{equation}\label{eq:known-LFL}
  \limsup_{n\to\infty}
  \frac{X_n}{(\log n)^{2/3}}
  =(2\beta)^{-2/3},
  \qquad
  \liminf_{n\to\infty}
  \frac{X_n}{(\log n)^{1/3}}
  =-\left(\frac8\beta\right)^{1/3}.
\end{equation}
The different logarithmic exponents reflect the
asymmetry of the upper and lower tails of the
Tracy-Widom distributions. This law is termed as
the \emph{law of fractional logarithm} (LFL) in ~\cite{PZ}.

In this paper, different from \cite{BCEHK-LFL},  we go beyond the GUE case and consider another direction of  extension:
keeping the Gaussian potential while allowing the
repulsion parameter $\beta$ to be any positive real
number. The natural setting is the Gaussian
$\beta$ corner process, which extends the classical
GOE and GUE minor processes to general $\beta>0$ \cite{Baryshnikov2001,JN,BorodinGorin2015,GS,Cuenca2021,GK,NV}. In the next subsection, we will introduce the model and the main results.

\subsection{Main results}\label{subsec:main}

Fix $\beta>0$. The (single level) Gaussian $\beta$ ensemble is a particle system defined on the chamber
$\cC_n:=\{x\in\mathbb R^n:x_1>\cdots>x_n\}$ with density
\begin{equation}\label{eq:rowlaw}
  \mu_{n,\beta}(\mathrm dx)
  =\frac{1}{Z_{n,\beta}}
   \exp\!\left(-\frac{\beta}{4}\sum_{i=1}^n x_i^2\right)
   \prod_{1\le i<j\le n}(x_i-x_j)^\beta\,\mathrm dx.
\end{equation}
It can be regarded as the jpdf of the eigenvalue of GOE, GUE and GSE in the case of $\beta=1,2,4$. For general $\beta>0$, it can also be regarded as the jpdf of eigenvalues of a tridiagonal random matrix model raised by \cite{DE}. However, the eigenvalues of the successive minors of the tridiagonal model does not give the corner process we need.  Instead, the Gaussian $\beta$ corner process we are going to consider  couples the particles in (\ref{eq:rowlaw}) from  different levels through interlacing as follows.  For
$x\in\cC_n$ and $y\in\cC_{n-1}$, write $y\prec x$ if
\[
  x_1>y_1>x_2>\cdots>y_{n-1}>x_n.
\]
With $\Delta(x):=\prod_{i<j}(x_i-x_j)$, the
Dixon-Anderson transition kernel is given by
\begin{equation}\label{eq:DA}
  \Lambda_{n,n-1}(x,\mathrm dy)
  =\frac{\Gamma(n\beta/2)}{\Gamma(\beta/2)^n}
   \frac{\Delta(y)}{\Delta(x)^{\beta-1}}
   \prod_{i=1}^{n-1}\prod_{j=1}^{n}
     |y_i-x_j|^{\beta/2-1}
   \mathbf 1_{\{y\prec x\}}\,\mathrm dy.
\end{equation}

\begin{defn}[Gaussian $\beta$ corner process]
\label{def:corners}
The Gaussian $\beta$ corner process is the infinite
interlacing array
$\{x_i^{(n)}:1\le i\le n,\ n\ge1\}$
whose first $N$ rows have joint law
\begin{equation}\label{eq:arraylaw}
  \mu_{N,\beta}(\mathrm dx^{(N)})
  \prod_{n=2}^{N}\Lambda_{n,n-1}
      (x^{(n)},\mathrm dx^{(n-1)}).
\end{equation}
\end{defn}

In particular, each row has law $\mu_{n,\beta}$,
while the dependence between successive rows is
specified by the Dixon-Anderson kernel. For
$\beta=1,2$, this process coincides with the
eigenvalue minor processes of the GOE and GUE,
respectively. We refer the reader to \cite{GS,Cuenca2021,GorinMarcus2020}
for further background on the Gaussian $\beta$ corner process.

For this corner process, we have the following LFL, as an extension of the results for the GOE/GUE corner process \cite{PZ,BBBK-LFL,BCEHK-LFL}.

\begin{thm}[Law of fractional logarithm]\label{thm:LFL}
For every fixed $\beta>0$, the Gaussian $\beta$ corner process
satisfies, almost surely,
\begin{equation}\label{eq:LFL}
  \limsup_{n\to\infty}\frac{X_n}{(\log n)^{2/3}}=(2\beta)^{-2/3},
  \qquad
  \liminf_{n\to\infty}\frac{X_n}{(\log n)^{1/3}}=-(8/\beta)^{1/3}.
\end{equation}
By symmetry, the same conclusions hold at the lower edge with
$X_n^{\mathrm{low}}:=n^{1/6}(-x_n^{(n)}-2\sqrt n)$ in place of $X_n$.
\end{thm}
Apart from the $\limsup$ and $\liminf$, our approach allows us to identify all the possible limit points under two different logarithmic scalings in (\ref{eq:LFL}), as \cite{BCEHK-LFL}. 
\begin{cor}[All limit points]\label{cor:limitpoints}
Almost surely, the sets of finite subsequential limits of
$X_n/(\log n)^{2/3}$ and $X_n/(\log n)^{1/3}$ are, respectively,
\begin{equation}\label{eq:limitpoints}
  [0,(2\beta)^{-2/3}]\qquad\text{and}\qquad[-(8/\beta)^{1/3},\infty).
\end{equation}
The same holds at the lower edge.
\end{cor}

For simplicity, we denote $\rR:=(2\beta)^{-2/3}$ and $\rL:=(8/\beta)^{1/3}$ in the remaining of the paper.

\normalcolor

\subsection{Proof ideas}
\label{subsec:proofideas}

On a high level, the overall proof strategy follows the framework of \cite{PZ,BBBK-LFL,BCEHK-LFL}, and is more directly related to 
~\cite{BCEHK-LFL}. The proof of LFL requires three inputs: a
tail probability estimates for the largest particles in a small deviation regime,
a short range correlation estimate for particles from neighboring levels, and long range
decorrelation estimates for particles from levels far apart. Similarly to \cite{PZ,BCEHK-LFL,BCEHK-Dec}, one expect that the correlation-decorrelation transition occurs when the level differences crosses $N^{2/3}$. This implies that  roughly
$N^{1/3}$ effectively independent observations are
among the levels in $[N,2N]$.
Together with the right and left tail probability in small deviation regime \cite{BBBK-Tails},
\[
  \log\P(X_n\ge t)\sim-\frac{2\beta}{3}t^{3/2},
  \qquad
  \log\P(X_n\le-t)\sim-\frac{\beta}{24}t^3,
\]
this explains the constants in
Theorem~\ref{thm:LFL}.
The correlation estimate controls the
oscillation between nearby levels, allowing
the upper bounds to be reduced to a finite
grid. Decorrelation provides approximate
independence on a sparse grid, which is
needed to establish the matching lower
bounds. We implement this framework in
Section~\ref{sec:mainproof}.

The main feature of our approach is that
both correlation and decorrelation can be
studied directly through the secular equation
dynamics of the Gaussian $\beta$ corner
process, without introducing an auxiliary
matrix evolution. We will simply term these dynamics as {\it secular dynamics} in the sequel. Actually, such kind of dynamic has been used in the martingale argument for the correlation part of LFL for Wigner minor process in \cite{BCEHK-LFL}. But the decorreltaion part in \cite{BCEHK-LFL} relies on a continuous dynamic, the Dyson Brownian motion.  

As the first ingredient, the tail probability estimates in small deviation regime, is already available thanks to the work \cite{BBBK-Tails}, together with Proposition \ref{prop:smallbeta} in the sequel, we will focus on introducing the strategy for proving the correlation-decorrelation transition. 

\medskip
\paragraph{Discrete secular dynamics.}
The Gaussian $\beta$ corner process admits
explicit secular equations in both directions.
For the \emph{upward transition} from level $n$ to
$n+1$, we use the representation in
\cite[Section~9]{NV}. Let
\begin{equation}\label{def.g-gamma}
  g^{(n)}\sim\mathcal N(0,2/\beta),
  \qquad
  \xi_i^{(n)}
  \sim\operatorname{Gamma}(\beta/2,2/\beta),
  \quad 1\le i\le n,
\end{equation}
be independent, where the Gamma distribution
is parametrized by shape and scale.
The particles at level $n+1$ are the ordered
roots of
\begin{equation}\label{eq:secularup}
  z-g^{(n)}
  -\sum_{i=1}^{n}
   \frac{\xi_i^{(n)}}{z-x_i^{(n)}}=0.
\end{equation}
Fresh independent innovations are used at
each step, where $x_i^{(n)}$'s are the level $n$ particles.

Conversely, the \emph{downward transition} from
level $n$ to $n-1$ admits the representation
in \cite[Equation~(4.9)]{Forrester2010}.
The particles at level $n-1$ are the
ordered roots of
\begin{equation}\label{eq:seculardown}
  \sum_{j=1}^{n}
  \frac{w_j^{(n)}}{z-x_j^{(n)}}=0,
\end{equation}
where
$w^{(n)}\sim\operatorname{Dirichlet}
(\beta/2,\ldots,\beta/2)$
is independent of the upper row.

\medskip

\paragraph{Long range decorrelation.}  In the case of the Wigner minor process \cite{BCEHK-LFL}, the authors used a dynamic approach for the decorrelation part, by running the same matrix-valued Brownian motion for a Wigner matrix of dimension $N+N^{2/3+\varepsilon}$ and also its minor of dimension $N$ (for the latter only the corresponding minor of the Brownian motion is used).  An eigenvector overlap estimate shows that the corresponding two correlated Dyson Brownian motions (DBMs) for two Wigner matrices of different dimensions get decoupled at the edge,  after a sufficiently long time, and thus the leading eigenvalues become independent. Then a Green function comparison leads to the decorrelation of the largest eigenvalues of the original matrices.  In practice, such a dynamic approach for  decorrelation  also  relies on a coupling for the two DBMs with two completely independent DBMs with independent initial conditions. Such a coupling idea for weakly correlated DBMs was first used in \cite{CES2023}, and we also refer to \cite{CES2021,BCH2024} etc for related discussion. Apart from application on the decorrelation of weakly correlated DBMs, the coupling idea was first introduced in \cite{BEYY} for the universality of local spectral statistics. We also refer to \cite{LandonYau2017,LandonSosoeYau2019} etc for further development.

Here, in our case, although we are also going to prove a decorrelation estimate, we nevertheless do not run any additional DBM starting from our initial particles.  Our approach to decorrelation is based on  the upward secular dynamics directly.
Starting from an independent Gaussian
$\beta$ ensemble row at level $N$, we construct
a second process using the same future innovations
as the original one. The restarted process is
independent of the history up to level $N$.
We show that the two evolutions become closer
than the natural edge fluctuation scale after
$N^{2/3+\varepsilon}$ steps. Consequently, we have the desired 
decorrelation estimate for the original process.
Conceptually, it can be regarded as a discrete version of the coupling approach in \cite{BEYY,LandonYau2017,LandonSosoeYau2019}.

The key observation is that subtracting the
two upward secular equations yields a 
{\it substochastic } evolution for their differences; see \eqref{eq:restart-difference}. Here by substochastic we mean that the $i$-sum of the coefficients $P_{ji}^{(n)}=C_{ji}^{(n)}/D_j^{(n)}$ is small than $1$, and thus $P=(P_{ji})$ is not a full stochastic transition matrix. We then track the dynamic of the maximal difference $\|\delta^{(n)}\|_{\infty}$. 
While substochasticity ensures that $\|\delta^{(n)}\|_{\infty}$ cannot increase, we will need a quantitative
contraction to show that it is sufficiently small after sufficiently many steps. According to eigenvalue rigidity, with high probability,  the large differences are from the spectral
edges. The main mechanism is that the upward transition can transfer
weights of the edge particle differences to those of the interior particles. The downward Dirichlet
representation ensures that such contractions of $\|\delta^{(n)}\|_{\infty}$
occur sufficiently often, i.e, for many $n$. The  accumulation of these contractions
gives the desired conclusion. 
\medskip

\paragraph{Short range correlation.} For the correlation estimate, our proof strategy is very similar to the counterpart of \cite{BCEHK-LFL}. More specifically, we again rely on 
the  secular dynamics. The upward Gamma  equation yields a
 one step estimates for
the increments of the largest particle.
After centering the innovations, the
remaining spectral terms are controlled
by eigenvalue rigidity. The accumulated
centered fluctuations can then be
estimated by martingale maximal
inequalities. Consequently, we can get a uniform bound
on the oscillation of $X_n$ over
intervals shorter than $N^{2/3}$.
See Section \ref{sec:correlation} for details.

\subsection{Notation.} Throughout the paper, 
We write $a\lesssim b$ if $a\le Cb$ for some
constant $C>0$, and use $a\gtrsim b$ and
$a\asymp b$ with the analogous meanings.
The constants $c,C>0$ may change from line to
line and may depend on $\beta$, the polynomial
decay exponent $D$ in probability estimates,
and fixed auxiliary parameters such as
$\varepsilon,\tau,\delta$.
Unless otherwise stated, they are independent
of the matrix size, level indices, and other
parameters varying in the stated ranges.

\normalcolor

\section{Proofs of the main results}\label{sec:mainproof}
Our proof follows the general framework of \cite{BCEHK-LFL}, relying on three ingredients: tail estimates, decorrelation, and correlation. We begin with the tail estimates for the rescaled largest particle $X_n$. The sharp decay rates of the right and left tails determine the candidate constants in the LFL. The following proposition states the needed tail estimates in the small deviation regime, which is essentially from \cite{BBBK-Tails}, except for the right tail upper bound in the regime $0<\beta<1$.

\begin{pro}[Tail probability estimates]\label{prop:tails}
Fix $K>0$ and $0<\varepsilon<1$. There are $T,n_0>0$ such that,
for $n\ge n_0$,
\begin{align}
 e^{-\frac{2\beta}{3}(1+\varepsilon)t^{3/2}}
 &\lesssim\P(X_n\ge t)
 \lesssim e^{-\frac{2\beta}{3}(1-\varepsilon)t^{3/2}},
 &&T\le t\le K(\log n)^{2/3},\label{eq:righttail}\\
 e^{-\frac\beta{24}(1+\varepsilon)t^3}
 &\lesssim\P(X_n\le-t)
 \lesssim e^{-\frac\beta{24}(1-\varepsilon)t^3},
 &&T\le t\le K(\log n)^{1/3}.\label{eq:lefttail}
\end{align}
\end{pro}

\begin{proof}
The right tail lower bound is Theorem~1.3(i) of
\cite{BBBK-Tails}.
The right tail upper bound for $\beta\ge2$ is their Theorem~1.1;
for $1\le\beta<2$ it is their Theorem~1.2(i).
Their Theorems~1.4(i) and~1.5(i) give, respectively, the left tail upper
and lower bounds for all $\beta>0$. The polynomial ranges of $t$
in those theorems contain the logarithmic ranges above for all sufficiently
large $n$. The remaining right tail upper bound, for $0<\beta<1$,
is proved in Proposition~\ref{prop:smallbeta} below.
\end{proof}

The tail estimates, Proposition \ref{prop:tails},  identify the candidate constants
in the LFL, but they are not sufficient to establish the
almost sure extremal behavior of the process. The main
difficulty lies in the dependence between different levels.
On one hand, strong correlations between nearby levels
prevent us from treating all indices as independent trials.
On the other hand, the tail estimates alone do
not guarantee that extreme excursions occur sufficiently
often.

Taking the $\limsup$ as an example, we shall  apply the first Borel-Cantelli lemmas to get a tight upper bound, but we cannot apply it to the full sequence directly due to the strong correlation of largest particles from neighboring levels. Hence, we shall first use it for a subsequence and then show the full sequence can be approximated by the subsequence. This would require us to prove a correlation estimate for the largest particles from neighboring levels. In contrast, for the tight lower bound, we shall apply the second Borel-Cantelli lemma to a sparser subsequence. To this end, we shall show that this sparser subsequence is almost independent, which is our decorrelation estimate. 

Recall $X_n$ from \eqref{eq:edgevariable}.
For $c>0$, we define the right and left tail events
\begin{equation}\label{eq:excursions}
 A_n^+(c):=\{X_n\ge c(\log n)^{2/3}\},\qquad
 A_n^-(c):=\{X_n\le-c(\log n)^{1/3}\}.
\end{equation}

We first state our decorrelation estimate. It shows that as long as the level difference is much larger than $N^{2/3}$, the tail events get  asymptotically independent.

\begin{pro}[Decorrelation estimate]\label{prop:decorrelation}
Fix $0<\delta<1/3$ and $K>1$. There is $\omega>0$ such that,
for $N\le m<n\le2N$ with $n-m\ge N^{2/3+\delta}$, and
$c_1,c_2\in[K^{-1},K]$, one has, for either choice of sign,
\begin{align}\label{eq:decorrelation}
 \P\bigl(A_m^\pm(c_1)\cap A_n^\pm(c_2)\bigr)
 &=\P\bigl(A_m^\pm(c_1)\bigr)
   \P\bigl(A_n^\pm(c_2)\bigr)
   \bigl(1+O(N^{-\omega}\log N)\bigr)
   +O(N^{-D})
\end{align}
for every $D>0$, uniformly in the stated parameters.
\end{pro}

For the correlation estimate, 
we will show that as long as the level difference is much smaller than
$N^{2/3}$, the fluctuation of the largest particles are essentially synchronized.

\begin{pro}[Correlation estimate]\label{prop:correlation}
For every $\varepsilon,D>0$, uniformly over
$1\le L\le N/2$ and $M=N-L$, we have
\begin{equation}\label{eq:correlation}
 \P\left(
 \max_{M\le a<b\le N}(X_a-X_b)>
 N^\varepsilon\left(
 \sqrt{\frac{L}{N^{2/3}}}
 +\frac{L}{N^{2/3}}
 \right)
 \right)
 \lesssim N^{-D}.
\end{equation}
In particular, for every fixed $0<\kappa<2/3$
and $D>0$,
\begin{equation}\label{eq:shortinterval}
 \P\left(
 \max_{0\le h_1<h_2\le\lfloor N^{2/3-\kappa}\rfloor}
 (X_{N+h_1}-X_{N+h_2})>1
 \right)
 \lesssim N^{-D}.
\end{equation}
\end{pro}

The proofs of Proposition \ref{prop:decorrelation} and \ref{prop:correlation} are deferred to Section
\ref{sec:decorrelation} and Sections~\ref{sec:correlation}. In the remaining part of this section,
we use the decorrelation estimate, correlation estimate and tail probability estimates to prove the main theorems.

\begin{proof}[Proof of Theorem \ref{thm:LFL}]
The proof is divided into two steps. 
We first establish the upper bound for the limsup
and the lower bound for the liminf by applying
the tail estimates along a sufficiently
dense subsequence and using the correlation estimate
to control the fluctuations between consecutive points.
For the reverse inequalities, it suffices to find
infinitely many excursions along a sparse subsequence,
where the decorrelation estimate provides approximate
independence and allows us to apply a general second Borel-Cantelli lemma.

\textit{Step 1: Upper bound for the limsup and lower bound for the liminf.}
Recall that $\rR:=(2\beta)^{-2/3}$ and $\rL:=(8/\beta)^{1/3}$.
We first prove that, almost surely,
\begin{equation}\label{eq:outergoal}
 \limsup_{n\to\infty}
 \frac{X_n}{(\log n)^{2/3}}\le\rR,
 \qquad
 \liminf_{n\to\infty}
 \frac{X_n}{(\log n)^{1/3}}\ge-\rL.
\end{equation}
We use a subsequence argument, followed by the
correlation estimate to control the fluctuations
between consecutive subsequence points.

Fix $C_R>\rR$ and $C_L>\rL$, and choose
$c_R\in(\rR,C_R)$ and $c_L\in(\rL,C_L)$.
Since
$
 2\beta\rR^{3/2}/3
 =\beta\rL^3/24=1/3,
$
we may choose $\varepsilon>0$ sufficiently small
so that
\[
 a_R:=\frac{2\beta}{3}(1-\varepsilon)c_R^{3/2}
 >\frac13,
 \qquad
 a_L:=\frac{\beta}{24}(1-\varepsilon)c_L^3
 >\frac13.
\]
Choose $p<3$ sufficiently close to $3$ such that
$
 pa_R>1
 $ and $pa_L>1$,
and consider the subsequence
\[
 n_k=\lfloor k^p\rfloor.
\]

By the tail estimates \eqref{eq:righttail}
and \eqref{eq:lefttail}, we have
\[
 \P\left(
 X_{n_k}\ge c_R(\log n_k)^{2/3}
 \right)
 \lesssim n_k^{-a_R}
 \lesssim k^{-pa_R},\quad \P\left(
 X_{n_k}\le-c_L(\log n_k)^{1/3}
 \right)
 \lesssim n_k^{-a_L}
 \lesssim k^{-pa_L}.
\]
Both bounds are summable in $k$. Hence,
by the Borel--Cantelli lemma, almost surely,
for all sufficiently large $k$,
\begin{equation}\label{eq:outer-subsequence}
 -c_L(\log n_k)^{1/3}
 <X_{n_k}
 <c_R(\log n_k)^{2/3}.
\end{equation}

It remains to control the fluctuations between $[n_k,n_{k+1}]$. Since $p<3$,
we can choose $\kappa>0$ such that
$
 1-p^{-1}<2/3-\kappa.
$
Observe that
\[
 n_{k+1}-n_k
 =O(k^{p-1})
 =O\left(n_k^{1-1/p}\right)
 =o\left(n_k^{2/3-\kappa}\right).
\]
Thus, Proposition~\ref{prop:correlation} gives,
for every $D>0$ and all sufficiently large $k$,
\[
 \P\left(
\max_{n_k\le a<b\le n_{k+1}}
(X_a-X_b)>1
\right)
\lesssim n_k^{-D}
\lesssim k^{-pD}.
\]
Taking $D$ large enough and applying Borel--Cantelli
once more, we obtain
\begin{equation*}
 \max_{n_k\le a<b\le n_{k+1}}
(X_a-X_b)\le1,\qquad\text{a.s.,}
\end{equation*}
for all sufficiently large $k$.
Combining this with \eqref{eq:outer-subsequence}, we obtain for every $n_k\le n\le n_{k+1}$ with $k$ large enough
\[
-c_L(\log n_k)^{1/3}-1
\le X_n
\le c_R(\log n_{k+1})^{2/3}+1,\qquad\text{a.s.}.
\]
Since $\log n_k\sim\log n_{k+1}$ and
$n_k\le n\le n_{k+1}$, it follows that
\[
 \frac{X_n}{(\log n)^{2/3}}
 \le c_R+o(1)<C_R,
 \qquad
 \frac{X_n}{(\log n)^{1/3}}
 \ge-c_L-o(1)>-C_L.
\]
Letting $C_R\downarrow\rR$ and
$C_L\downarrow\rL$ along countable sequences
completes the proof of \eqref{eq:outergoal}.

\textit{Step 2:Lower bound for the limsup and upper bound for the liminf.}
We now prove the reverse inequalities by constructing a subsequence satisfying the desired bound, following the classical idea of the second
Borel--Cantelli lemma.
Fix $0<c_R<\rR$ and $0<c_L<\rL$. Choose
$\varepsilon>0$ sufficiently small so that
\[
 a_R:=\frac{2\beta}{3}(1+\varepsilon)c_R^{3/2}
 <\frac13,
 \qquad
 a_L:=\frac{\beta}{24}(1+\varepsilon)c_L^3
 <\frac13.
\]
We may then choose $p>3$ sufficiently close to $3$
such that $pa_R<1$ and $pa_L<1$.
Consider the subsequence
\[
 n_k=\lfloor k^p\rfloor.
\]

For simplicity, we write $c_+=c_R, c_-=c_L$ and $a_+=a_R, a_-=a_L$  either sign $\sigma\in\{+,-\}$ for either sign $\sigma\in\{+,-\}$.
We group
the subsequence into blocks on which the events
are approximately independent. For each sufficiently
large integer $M$, set
\[
I_M:=\{k:n_M\le n_k<2n_M\}.
\]
It follows from $n_k=\lfloor k^p\rfloor$ that
$|I_M|\asymp n_M^{1/p}$. As $p>3$, we may choose $\delta>0$ such that
$1-p^{-1}>2/3+\delta$. Therefore, any two distinct
indices $k,\ell\in I_M$ satisfy
\[
|n_k-n_\ell|\gtrsim n_M^{1-1/p}\gtrsim n^{2/3+\delta}_M.
\]
Thus, for all
sufficiently large $M$, Proposition~\ref{prop:decorrelation} is applicable.

On the other hand, the tail estimates
\eqref{eq:righttail} and \eqref{eq:lefttail}
give
\begin{equation}\label{eq:manyexcursions}
\sum_{k\in I_M}\P(A_{n_k}^\sigma(c_\sigma))
\gtrsim |I_M|n_M^{-a_\sigma}
\asymp n_M^{1/p-a_\sigma}
\xrightarrow{M\to\infty}\infty,
\end{equation}
where we used $pa_\sigma<1$.

By the Chung-Erd\H{o}s inequality
\cite[Lemma, p.~180, (4)]{ChungErdos1952},
we have
\[
\P\left(
\bigcup_{k\in I_M}A_{n_k}^\sigma(c_\sigma)
\right)
\ge
\frac{
\left(\sum_{k\in I_M}
\P(A_{n_k}^\sigma(c_\sigma))\right)^2
}{
\displaystyle\sum_{k,\ell\in I_M}
\P\left(
A_{n_k}^\sigma(c_\sigma)
\cap A_{n_\ell}^\sigma(c_\sigma)
\right)
}.
\]

For simplicity, write
$S_M:=\sum_{k\in I_M}\P(A_{n_k}^\sigma(c_\sigma))$.
We compare the denominator with $S_M^2$
by separating the diagonal terms and the
correlations between distinct events:
\begin{align*}
\sum_{k,\ell\in I_M}\P(A_{n_k}^\sigma(c_\sigma)\cap A_{n_l}^\sigma(c_\sigma))
&=S_M^2+\sum_{k\in I_M}\P(A_{n_k}^\sigma(c_\sigma))(1-\P(A_{n_k}^\sigma(c_\sigma)))\\
&+\sum_{\substack{k,\ell\in I_M\\k\ne\ell}}
\bigl(\P(A_{n_k}^\sigma(c_\sigma)\cap A_{n_l}^\sigma(c_\sigma))-\P(A_{n_k}^\sigma(c_\sigma))\P(A_{n_l}^\sigma(c_\sigma))\bigr)\\
&\ge \bigl(1-CM^{-\omega}\log M\bigr)S_M^2
      -C|I_M|^2M^{-D},
\end{align*}
where we used the nonnegativity of the diagonal contribution and Proposition \ref{prop:decorrelation}.
Putting the above inequality back into the Chung-Erd\H{o}s inequality and using
the elementary inequality $(1+x)^{-1}\ge 1-x$ for $x\ge0$ gives
\begin{equation}\label{eq:secondmoment}
\begin{aligned}
\P\left(\bigcup_{k\in I_M}A_k\right)
\ge
1-O\left(
N_M^{a_\sigma-1/p}
+N_M^{-\omega}\log N_M
+N_M^{2a_\sigma-D}
\right)
\xrightarrow{M\to\infty}1,
\end{aligned}
\end{equation}
where we used \eqref{eq:manyexcursions}
and $|I_M|\asymp N_M^{1/p}$.

For every fixed $K$, we have
$I_M\subset\{k:k\ge K\}$ for all sufficiently
large $M$. Hence \eqref{eq:secondmoment} implies
\[
\P\Big(
\bigcup_{k\ge K}A_{n_k}^\sigma(c_\sigma)
\Big)=1.
\]
Taking the intersection over $K\in\mathbb N$,
we conclude that
$
\P\left(
A_{n_k}^\sigma(c_\sigma)\ \mathrm{i.o.}
\right)=1.$

Consequently, almost surely, we have
\[
 \limsup_{n\to\infty}
 \frac{X_n}{(\log n)^{2/3}}\ge c_R,
 \qquad
 \liminf_{n\to\infty}
 \frac{X_n}{(\log n)^{1/3}}\le-c_L.
\]
Sending  $c_R$ and $c_L$ to $\rR$ and
$\rL$ respectively, we obtain the desired bounds.
Together with \eqref{eq:outergoal}, this
completes the proof of \eqref{eq:LFL}.

\end{proof}

The LFL determines the asymptotic upper and lower
envelopes of $X_n$. We now show that every value
between these envelopes can be attained as a
subsequential limit after normalization. The key
observation is that the correlation estimate makes
the increments of the normalized sequence tend to
zero, so that it cannot skip any intermediate value.
The lower edge results then follow by symmetry.

\begin{proof}[Proof of Corollary~\ref{cor:limitpoints}]
We first consider the upper edge. Since $D$ is large enough, the probabilities in
\eqref{eq:shortinterval} is summable. Hence
Borel-Cantelli lemma implies
\begin{equation}\label{eq:smalljumps}
\P\left(
X_{n+1}-X_n\ge -1
\text{ for all sufficiently large }n
\right)=1.
\end{equation}
Moreover, the LFL implies that
$|X_n|=O((\log n)^{2/3})$ almost surely.
Consequently, for $s=1/3$ or $2/3$, we have
\[
\begin{aligned}
\frac{X_{n+1}}{(\log(n+1))^s}
-\frac{X_n}{(\log n)^s}
&=
\frac{X_{n+1}-X_n}{(\log(n+1))^s}
+X_n\left(
\frac{1}{(\log(n+1))^s}
-\frac{1}{(\log n)^s}
\right)
\\
&\ge
-\frac{1}{(\log(n+1))^s}
-\frac{s|X_n|}{n(\log n)^{s+1}}
\xrightarrow{\mathrm{a.s.}}0,
\end{aligned}
\]
which directly gives
\begin{align}
\liminf_{n\to\infty}
\left(
\frac{X_{n+1}}{(\log(n+1))^s}
-\frac{X_n}{(\log n)^s}
\right)\ge0
\qquad\text{a.s.} \label{111001}
\end{align}

Next, by  the  LFL in \eqref{eq:LFL}, we also have 
\[
 \liminf_{n\to\infty}
 \frac{X_n}{(\log n)^{2/3}}=0, \qquad\limsup_{n\to\infty}
 \frac{X_n}{(\log n)^{1/3}}=+\infty.
\]
From the above two facts, we can actually conclude that 
 every  value between the liminf and limsup is a  limit point. This is because that the sequence will go from the limsup down to liminf for infinitely many often, but (\ref{111001}) on the liminf of the increments shows that one cannot go down too fast at each step to skip any possible limit point.  This completes the proof for the upper edge. The lower edge case follows by symmetry. 
\end{proof}

\section{Proof of the decorrelation estimate} \label{sec:decorrelation}

In this section, we prove Proposition \ref{prop:decorrelation}. As mentioned in Section \ref{subsec:proofideas},  our strategy is to couple the original process with a process started independently at  level $N$, using the same innovations afterwards in the upward transitions. 
We then show that after more than $N^{2/3}$  steps, the two processes become sufficiently close. Since the auxiliary process is independent of the history before level $N$, this closeness yields decorrelation between level $N$ and level $N+N^{2/3+\varepsilon}$.

More precisely, let $\cF_N$ denote the history of the original process through level $N$. We choose $\widetilde x^{(N)}\sim\mu_{N,\beta}$ independently of $\cF_N$ and all future innovations, and construct $\widetilde x^{(n)}$ for $n>N$ by iterating the upward transition $n-N$ times, using the same innovations as the original process at each step. Specifically, for each $n\ge N$, let
\[
g^{(n)}\sim N(0,2/\beta),\qquad
\xi_i^{(n)}\sim\GammaLaw(\beta/2,2/\beta),
\quad 1\le i\le n,
\]
be independent. The next rows $x^{(n+1)}$ and $\widetilde x^{(n+1)}$ are then obtained as the ordered roots of
\begin{equation}\label{eq.secular equation-both}
z-g^{(n)}-\sum_{i=1}^n
\frac{\xi_i^{(n)}}{z-x_i^{(n)}}=0,
\qquad
z-g^{(n)}-\sum_{i=1}^n
\frac{\xi_i^{(n)}}{z-\widetilde x_i^{(n)}}=0,
\end{equation} respectively.

We can think of this as coupling two Markov chains starting from different initial configurations by driving them with the same innovations. After sufficiently many upward steps, the two chains largely forget their initial configurations and become very close to each other.
The following proposition establishes this quantitative coupling estimate.

\begin{pro}[Independent restart]\label{prop:restart}
Fix $0<\varepsilon<1/3$. There exists $\omega>0$ such that,
uniformly for $N^{2/3+\varepsilon}\le k\le N$,
\begin{equation}\label{eq:restart}
 \P\left(
 \max_i|x_i^{(N+k)}-\widetilde x_i^{(N+k)}|
 >N^{-1/6-\omega}
 \right)\lesssim N^{-D}
\end{equation}
for every $D>0$ and all sufficiently large $N$.
Moreover, $\widetilde x^{(N+k)}$ has law $\mu_{N+k,\beta}$
and is independent of $\cF_N$.
\end{pro}

With the above proposition in hand, we can now directly prove Proposition \ref{prop:decorrelation}.

\begin{proof}[Proof of Proposition~\ref{prop:decorrelation},
assuming Proposition~\ref{prop:restart}]
Restart the process at level $m$ and evolve both copies
with the same innovations up to level $n$.
Since $n-m\ge N^{2/3+\delta}$ and $m\asymp N$,
Proposition~\ref{prop:restart}, applied with
$\varepsilon=\delta/2$, gives
\begin{equation}\label{eq:edgecoupling}
\P\bigl(|X_n-\widetilde X_n|>N^{-\omega_0}\bigr)
\lesssim N^{-D},
\end{equation}
where
$\widetilde X_n=n^{1/6}(\widetilde x_1^{(n)}-2\sqrt n)$
has the same law as $X_n$ and is independent of $\cF_m$.

For $\sigma\in\{+,-\}$, write
\[
Q_n^+(t):=\P(X_n\ge t),\qquad
Q_n^-(t):=\P(X_n\le-t),
\]
and set $s_+=2/3$, $s_-=1/3$.
Let $t=c_2(\log n)^{s_\sigma}$ .
Independence and \eqref{eq:edgecoupling} yields
\begin{align*}
\P(A_m^\sigma(c_1))Q_n^\sigma(t+N^{-\omega_0})-CN^{-D}
\le
\P(A_m^\sigma(c_1)\cap A_n^\sigma(c_2))\le
\P(A_m^\sigma(c_1))Q_n^\sigma(t-N^{-\omega_0})+CN^{-D}.
\end{align*}

Since the row density \eqref{eq:rowlaw} is log-concave
for every $\beta>0$, Proposition~\ref{prop:tails},
together with the argument of \cite[Lemma~3.4]{BCEHK-LFL} can show that $Q_n^\sigma(t\pm N^{-\omega_0})$ can be approximated by $Q_n^{\sigma}(t)$ up to negligible error.  Consequently, we will have
\[
\P(A_m^\sigma(c_1)\cap A_n^\sigma(c_2))
=
\P(A_m^\sigma(c_1))\P(A_n^\sigma(c_2))
\bigl(1+O(N^{-\omega_1})\bigr),
\]
for some $\omega_1>0$. This proves \eqref{eq:decorrelation} for both tails.

\end{proof}

We now prove Proposition~\ref{prop:restart} by
tracking the maximal coordinate difference
between the two coupled processes. Recall that
$x^{(n)}$ and $\widetilde x^{(n)}$ evolve according
to the upward secular equations
\eqref{eq.secular equation-both}, driven by the
same innovations $g^{(n)}$ and
$\{\xi_i^{(n)}\}_{i=1}^n$.

Set
\[
  \delta_i^{(n)}
  :=x_i^{(n)}-\widetilde x_i^{(n)},
  \qquad N\le n\le 2N.
\]
Taking the difference of two secular equations of two consecutive levels gives
\begin{equation}\label{eq:restart-difference}
  \delta_j^{(n+1)}
  =\sum_{i=1}^n
   \frac{C_{ji}^{(n)}}{D_j^{(n)}}
   \delta_i^{(n)},\qquad \mbox{where}\quad C_{ji}^{(n)}
  :=\frac{\xi_i^{(n)}}{
  (x_j^{(n+1)}-x_i^{(n)})
  (\widetilde x_j^{(n+1)}
       -\widetilde x_i^{(n)})},
  \qquad
  D_j^{(n)}:=1+\sum_{i=1}^n C_{ji}^{(n)}.
\end{equation}
By interlacing, the two denominator factors
in $C_{ji}^{(n)}$ have the same sign, so
$C_{ji}^{(n)}>0$. Consequently, with
$P_{ji}^{(n)}:=C_{ji}^{(n)}/D_j^{(n)}$, we have
\[
  \delta_j^{(n+1)}
  =\sum_{i=1}^nP_{ji}^{(n)}\delta_i^{(n)},
  \qquad
  \sum_{i=1}^nP_{ji}^{(n)}
  =1-\frac1{D_j^{(n)}}<1.
\]
Thus the differences evolve through a
positive substochastic transition kernel. In
particular, this implies 
\begin{align}
  \|\delta^{(n+1)}\|_\infty\le \|\delta^{(n)}\|_\infty. \label{111002}
\end{align}
To obtain quantitative decay, we introduce the notions of 
{\it regular} and {\it effective} steps, as specified in
Definition~\ref{def:restart-steps}.
Regularity combines eigenvalue rigidity with
concentration of the Gamma innovations, while
effectiveness controls the normalization
factors $D_j^{(n)}$ near the spectral edges. If a step does not satisfy these properties, we may use the crude nonquantitative decay (\ref{111002}). 
The argument is based on two propositions.
Proposition~\ref{prop:restart-step} shows that
at every regular and effective step, either
$\|\delta^{(n)}\|_\infty$ has already reached the target scale
or a quantitative contraction occurs.
Proposition~\ref{prop:restart-count}, using
the downward Dirichlet representation, shows
that at least half of the steps are effective
with overwhelming probability. Since regularity
holds throughout the relevant interval with
overwhelming probability and $\|\delta^{(n)}\|_\infty$ cannot
increase at any step, these contractions
accumulate to give the desired estimate.

To prove the preceding two propositions, we introduce
an auxiliary parameter $K\in\mathbb N_+$ with
$1\le K\le N/8$, which will be chosen as a small
power of $N$ in the proof of Proposition~\ref{prop:restart}.
We also fix the following notation. Let
$\{q_i^{(r)}\}_i$ denote the
quantiles of the semicircle law, defined by
\begin{equation}\label{def.quantile}
q_i^{(r)}=\sqrt r\,\gamma_i^{(r)},\qquad\mbox{where}\quad
\int_{\gamma_i^{(r)}}^2
\frac{\sqrt{4-x^2}}{2\pi}\,\dd x
=\frac{i-1/2}{r},
\end{equation}
and set $\widehat i_r=\min\{i,r+1-i\}$.
Since our analysis focuses primarily on the spectral edges, we introduce the following index sets:
\[
E_r:=\{j:\widehat j_r\le K\},\qquad
J^+:=\{2K,\ldots,3K-1\},\qquad
J_n^-:=\{n+1-i:i\in J^+\}.
\]

We next introduce two conditions on the one step evolution that will be used to establish quantitative contraction.
\begin{defn}[Regular and effective steps]\label{def:restart-steps}
Let $N\le n<2N$. The transition from level $n$ to level $n+1$ is
\emph{regular} if
\begin{itemize}
\item[(i)] for $r\in\{n,n+1\}$ and $1\le i\le r$,
\begin{equation}\label{eq:restart-rigidity}
 |x_i^{(r)}-q_i^{(r)}|+|\widetilde x_i^{(r)}-q_i^{(r)}|
 \le K^{1/6}N^{-1/6}\,\widehat i_r^{-1/3};
\end{equation}
\item[(ii)] $\sum_{i\in J^+}\xi_i^{(n)}\ge K/2$ and
 $\sum_{i\in J_n^-}\xi_i^{(n)}\ge K/2$.
\end{itemize}
It is \emph{effective} if
\begin{equation}\label{eq:restart-effective}
 D_j^{(n)}\le NK^{4/\beta}
 \qquad\text{for every }j\in E_{n+1}.
\end{equation}
\end{defn}
Condition~(i) is a rigidity assumption, whose parameter will be specified at the end of the proof. It directly yields
$|\delta_i^{(r)}|\le K^{1/6}N^{-1/6},\widehat i_r^{-1/3}.$ Condition~(ii) ensures that the reference blocks carry
a proportional amount of innovation weight. Effectiveness concerns the
normalization of the weights $P_{ji}^{(n)}=C_{ji}^{(n)}/D_j^{(n)}$ in
\eqref{eq:restart-difference}: it prevents the contribution of the
reference blocks to the numerators from being outweighed by a large
denominator.
Under these conditions, we obtain the following one-step contraction estimate.

\begin{pro}[One step contraction]\label{prop:restart-step}
There is $c>0$ such that the following holds for $N\le n<2N$. If the
transition from level $n$ to level $n+1$ is regular and effective, then
\begin{equation}\label{eq:restart-step}
 \|\delta^{(n+1)}\|_\infty
 \le
 \max\Bigl\{
 2(NK)^{-1/6},\,
 \bigl(1-cN^{-2/3}K^{-1/3-4/\beta}\bigr)\|\delta^{(n)}\|_\infty
 \Bigr\}.
\end{equation}
\end{pro}

\begin{proof}
We may assume that $\|\delta^{(n+1)}\|_\infty>2(NK)^{-1/6}$,
since otherwise the claim is immediate.
Choose $j$ such that $|\delta_j^{(n+1)}|=\|\delta^{(n+1)}\|_\infty$.
It follows from \eqref{eq:restart-rigidity} that
\[
2(NK)^{-1/6}<\|\delta^{(n+1)}\|_\infty
\le K^{1/6}N^{-1/6}\widehat j_{n+1}^{-1/3},
\]
which implies $\widehat j_{n+1}<K/8$.
By symmetry, it suffices to consider $j<K/8$;
the opposite edge is treated in the same way
with $J_n^-$ in place of $J^+$.

Set
\[
I_n:=\{i:|\delta_i^{(n)}|\le \|\delta^{(n)}\|_\infty/2\}.
\]
Since the weights $P_{ji}^{(n)}$ are nonnegative
and have total mass at most one, we can obtain
\begin{equation}\label{eq:restart-split}
\begin{aligned}
\|\delta^{(n+1)}\|_\infty=|\delta^{(n+1)}_j|=&\left|\sum_{i\notin I_n}P_{ji}^{(n)}\delta^{(n)}_j+\sum_{i\in I_n}P_{ji}^{(n)}\delta^{(n)}_j\right|\\
\le&\|\delta^{(n)}\|_\infty\sum_{i\notin I_n}P_{ji}^{(n)}
  +\frac{\|\delta^{(n)}\|_\infty}{2}\sum_{i\in I_n}P_{ji}^{(n)}
\le \Big(1-\frac12
  \sum_{i\in I_n}P_{ji}^{(n)}\Big)\|\delta^{(n)}\|_\infty.
\end{aligned}
\end{equation}
We next bound the sum $\sum_{i\in I_n}P_{ji}^{(n)}$ from below by restricting it to $J^+$. To this end, we first show that $J^+\subset I_n$ and then estimate the corresponding weights.

It follows from $i\in J^+$ that
$2K\le i\le3K-1\le n/2$. Since the transition from level $n$ to $n+1$ is regular, 
condition (i) of Definition \ref{def:restart-steps} gives
\[
|\delta_i^{(n)}|
\le K^{1/6}N^{-1/6}(2K)^{-1/3}
<(NK)^{-1/6}
<\frac{\|\delta^{(n+1)}\|_\infty}2
\le\frac{\|\delta^{(n)}\|_\infty}2,
\]
and hence $J^+\subset I_n$.

Recall that the semicircle quantiles satisfy
\begin{equation}\label{eq.semi-location}
0\le 2\sqrt r-q_\ell^{(r)}
\lesssim r^{-1/6}\ell^{2/3},
\qquad 1\le\ell\le r/2.
\end{equation}
For $j<K/8$ and $i\in J^+$, we decompose
\[
\begin{aligned}
x_j^{(n+1)}-x_i^{(n)}
=&(x_j^{(n+1)}-q_j^{(n+1)})
 +(q_j^{(n+1)}-q_i^{(n)})+(q_i^{(n)}-x_i^{(n)})\\
 =&(x_j^{(n+1)}-q_j^{(n+1)})
 +(q_i^{(n)}-x_i^{(n)})+(q_j^{(n+1)}-2\sqrt{n+1})\\
 &+(2\sqrt{n+1}-2\sqrt n)
 +(2\sqrt n-q_i^{(n)})
\end{aligned}
\]
and similarly for $\widetilde x^{(n+1)}_j-\widetilde x^{(n)}_i$. 
Thus, by condition (i) of Definition \ref{def:restart-steps} and the
basic fact \eqref{eq.semi-location}, it can be derived that
\[
|x_j^{(n+1)}-x_i^{(n)}|
+|\widetilde x_j^{(n+1)}-\widetilde x_i^{(n)}|
\le CN^{-1/6}K^{2/3}.
\]
Consequently,
$C_{ji}^{(n)}\gtrsim
N^{1/3}K^{-4/3}\xi_i^{(n)}$,
and condition~(ii) of
Definition~\ref{def:restart-steps} yields
\begin{equation}\label{eq:restart-numerator}
\sum_{i\in J^+}C_{ji}^{(n)}
\gtrsim N^{1/3}K^{-1/3}.
\end{equation}

Finally, since $j\in E_{n+1}$,
effectiveness in Definition \ref{def:restart-steps} gives $D_j^{(n)}\le NK^{4/\beta}$.
Together with $J^+\subset I_n$ and
\eqref{eq:restart-numerator}, this implies
\[
\sum_{i\in I_n}P_{ji}^{(n)}
\ge \frac{1}{D_j^{(n)}}
    \sum_{i\in J^+}C_{ji}^{(n)}
\gtrsim N^{-2/3}K^{-1/3-4/\beta}.
\]
Substituting this into \eqref{eq:restart-split}
completes the proof.
\end{proof}

We next show that effective steps occur sufficiently often to ensure a cumulative contraction.

\begin{pro}[Effective steps]\label{prop:restart-count}
There is $C>0$ such that, for all sufficiently large $N$,
$1\le K\le N/8$, and $1\le L\le N$,
\begin{equation}\label{eq:restart-count}
\P\left(
\#\{N\le n<N+L:\text{the transition }n\to n+1
\text{ is effective}\}<\frac L2
\right)
\le 2\left(\frac CK\right)^{L/4}.
\end{equation}
\end{pro}

\begin{proof}
We first reduce the problem to the weights of the
individual processes. For $N\le n<N+L$, define
\begin{equation}\label{eq:restart-r}
r_j^{(n)}
:=\left(1+\sum_{i=1}^n
\frac{\xi_i^{(n)}}
{(x_j^{(n+1)}-x_i^{(n)})^2}\right)^{-1},
\qquad 1\le j\le n+1,
\end{equation}
and define $\widetilde r_j^{(n)}$ analogously
for the independent restarted process. By the Cauchy-Schwarz
inequality,
\[
\begin{aligned}
D_j^{(n)}
&=1+\sum_{i=1}^n
\frac{\xi_i^{(n)}}
{(x_j^{(n+1)}-x_i^{(n)})
(\widetilde x_j^{(n+1)}-\widetilde x_i^{(n)})}\\
&\le
\left(1+\sum_{i=1}^n
\frac{\xi_i^{(n)}}
{(x_j^{(n+1)}-x_i^{(n)})^2}\right)^{1/2}
\left(1+\sum_{i=1}^n
\frac{\xi_i^{(n)}}
{(\widetilde x_j^{(n+1)}
-\widetilde x_i^{(n)})^2}\right)^{1/2}=\bigl(r_j^{(n)}\widetilde r_j^{(n)}
\bigr)^{-1/2}.
\end{aligned}
\]
We introduce the events
\[
\begin{aligned}
B_n:=
\left\{\min_{j\in E_{n+1}}r_j^{(n)}
<N^{-1}K^{-4/\beta}\right\},\qquad
\widetilde B_n:=
\left\{\min_{j\in E_{n+1}}\widetilde r_j^{(n)}
<N^{-1}K^{-4/\beta}\right\}.
\end{aligned}
\]
On $B_n^c\cap\widetilde B_n^c$, we have
$$\min_{j\in E_{n+1}}\{r_j^{(n)},\widetilde r_j^{(n)}\}
\ge N^{-1}K^{-4/\beta}.$$ The preceding
Cauchy-Schwarz bound then gives
$D_j^{(n)}\le NK^{4/\beta}$, so the step
is effective. Consequently, every non-effective
step belongs to $B_n\cup\widetilde B_n$. Thus, it suffices to prove
\begin{equation}\label{eq:restart-count-one}
\P\left(
\sum_{n=N}^{N+L-1}\one_{B_n}\ge\frac L4
\right)
\le\left(\frac CK\right)^{L/4},
\end{equation}
together with the analogous estimate for
$\widetilde B_n$. We establish the former,
since both processes have the same marginal law.

The key observation is that the vectors $r^{(n)}$
are precisely the weights in the downward
Dirichlet representation. To see this, let
$p_n(w)=\prod_{i=1}^n(w-x_i^{(n)})$.
By the upward secular equation
\eqref{eq.secular equation-both}, multiplying its
left hand side by $p_n(w)$ yields a monic polynomial
of degree $n+1$, whose zeros are precisely
$x_1^{(n+1)},\ldots,x_{n+1}^{(n+1)}$.
This polynomial must therefore equal $p_{n+1}(w)$,
giving
\begin{equation}\label{eq:restart-charpoly}
\frac{p_{n+1}(w)}{p_n(w)}
=w-g^{(n)}
-\sum_{i=1}^n
\frac{\xi_i^{(n)}}{w-x_i^{(n)}}.
\end{equation}
Differentiating \eqref{eq:restart-charpoly}
at $w=x_j^{(n+1)}$ and using
$p_{n+1}(x_j^{(n+1)})=0$, we obtain
\[
\frac{p_{n+1}'(x_j^{(n+1)})}
{p_n(x_j^{(n+1)})}
=
1+\sum_{i=1}^n
\frac{\xi_i^{(n)}}
{(x_j^{(n+1)}-x_i^{(n)})^2}
=\frac1{r_j^{(n)}}.
\]
Since $x_j^{(n+1)}$ is a simple zero of
$p_{n+1}$, we have
\begin{align*}
\operatorname{Res}_{w=x_j^{(n+1)}}
\frac{p_n(w)}{p_{n+1}(w)}
&=\lim_{w\to x_j^{(n+1)}}
(w-x_j^{(n+1)})
\frac{p_n(w)}{p_{n+1}(w)}=\frac{p_n(x_j^{(n+1)})}
{p_{n+1}'(x_j^{(n+1)})}
=r_j^{(n)}.
\end{align*}
As $p_n/p_{n+1}$ has only simple poles
and vanishes at infinity, its partial fraction
expansion is given by
\begin{equation}\label{eq:restart-partialfractions}
\frac{p_n(w)}{p_{n+1}(w)}
=\sum_{j=1}^{n+1}
\frac{\displaystyle
\operatorname{Res}_{z=x_j^{(n+1)}}
\frac{p_n(z)}{p_{n+1}(z)}}
{w-x_j^{(n+1)}}=\sum_{j=1}^{n+1}
\frac{r_j^{(n)}}{w-x_j^{(n+1)}}.
\end{equation}
By the downward Markov structure
\eqref{eq:arraylaw}, conditional on the top row,
the lower rows can be constructed successively
using independent Dirichlet weight vectors.
Since these weights are uniquely determined
by adjacent rows through
\eqref{eq:restart-partialfractions}, the vectors
$r^{(N)},\ldots,r^{(N+L-1)}$ have the same
joint law as these independent Dirichlet vectors.
Consequently, they are independent, satisfying
\[
r^{(n)}\sim\Dir(\beta/2,\ldots,\beta/2).
\]
Since each coordinate of a Dirichlet vector
has a Beta distribution with parameters given
by its own parameter and the sum of the remaining
parameters, we have
$
r_j^{(n)}\sim
\operatorname{Beta}\left(\beta/2, n\beta/2\right).
$

We now estimate the probability of $B_n$.
Notice that the density of  $r_j^{(n)}$ is given by
\[
f_n(s)=
\frac{\Gamma((n+1)\beta/2)}
{\Gamma(\beta/2)\Gamma(n\beta/2)}
s^{\beta/2-1}(1-s)^{n\beta/2-1},
\qquad 0<s<1.
\]
Using the standard estimate
$\Gamma(x+\beta/2)/\Gamma(x)\lesssim x^{\beta/2}$
for $x\ge1$, we obtain, uniformly for
$N\le n<2N$ and $0<t\le N^{-1}$,
\[
\begin{aligned}
\P(r_j^{(n)}<t)
=\int_0^t f_n(s)\,\dd s\lesssim n^{\beta/2}
\int_0^t s^{\beta/2-1}\,\dd s
\lesssim (nt)^{\beta/2},
\end{aligned}
\]
where we used $(1-s)^{n\beta/2-1}\le1$
for sufficiently large $N$.
Taking $t=N^{-1}K^{-4/\beta}$ gives
\[
\P\left(
r_j^{(n)}<N^{-1}K^{-4/\beta}
\right)\lesssim K^{-2}.
\]
Since $|E_{n+1}|\le2K$, a union bound yields
$
\P(B_n)\lesssim K^{-1}.
$
Moreover, the events $B_n$ are independent,
as each depends only on the Dirichlet vector
$r^{(n)}$.

It remains to count the number of bad steps. By independence
and $
\P(B_n)\lesssim K^{-1}.
$,
\[
\begin{aligned}
\P\left(
\sum_{n=N}^{N+L-1}\one_{B_n}\ge\frac L4
\right)
\le\binom{L}{L/4}\left(\frac CK\right)^{L/4}\le 2^L\left(\frac CK\right)^{L/4}
\le\left(\frac{C'}K\right)^{L/4}.
\end{aligned}
\]
Here we may assume that $K$ is sufficiently
large, since otherwise the desired estimate
is trivial after enlarging the constant.
This proves \eqref{eq:restart-count-one}.
The same argument applies to $\widetilde B_n$,
and \eqref{eq:restart-count} follows by a
union bound.
\end{proof}

\begin{proof}[Proof of Proposition~\ref{prop:restart}]
Fix $0<\varepsilon<1/3$ and choose
\[
\kappa:=\frac{\beta\varepsilon}{\beta+12},
\qquad
K:=\lfloor N^\kappa\rfloor,
\qquad
k_0:=\lceil N^{2/3+\varepsilon}\rceil,
\qquad
\omega:=\frac{\kappa}{12}.
\]
With the above parameter choice, for sufficiently large $N$, we have 
\begin{equation}\label{eq:restart-exponent}
1\le K\le N/8,\qquad k_0\le N,\qquad
K^{1/3+4/\beta}
\le N^{\kappa(\beta+12)/(3\beta)}
=N^{\varepsilon/3}.
\end{equation}

Let
\[
\mathcal G:=
\left\{
\begin{array}{l}
\text{all transitions }n\to n+1,\ N\le n<N+k_0,
\text{ are regular,}\\
\text{at least }k_0/2\text{ of them are effective}
\end{array}
\right\}.
\]

We first show that $\mathcal G$ holds with overwhelming
probability. Since both processes $\{x^{(n)}\}$ and
$\{\widetilde x^{(n)}\}$ have the Gaussian
$\beta$ ensemble law $\mu_{r,\beta}$ at every level $r$,
the rigidity estimate \cite[Theorem~2.4]{BEY} yields,
for every $D>0$,
\[
\P\left(
|x_i^{(r)}-q_i^{(r)}|
>r^{-1/6+\kappa/12}\widehat i_r^{-1/3}
\text{ for some }1\le i\le r
\right)
\lesssim r^{-D},
\]
and the same bound holds for the restarted process
$\{\widetilde x^{(n)}\}$.
Since $K\ge N^\kappa/2$, taking a union bound
over both processes and all levels
$N\le r\le N+k_0$ shows that
\eqref{eq:restart-rigidity} holds simultaneously
at these levels with probability at least
$1-O(N^{-D})$. In particular, on this event, we have
\begin{equation}\label{eq.1}
\|\delta^{(n)}\|_\infty
\le K^{1/6}N^{-1/6}
\le N^{-1/6+\kappa/6}.
\end{equation}

For condition~(ii) in Definition~\ref{def:restart-steps},
each block contains $K$ indices, and the
corresponding innovation weights are independent
Gamma variables with mean one. A standard
concentration estimate gives
\[
\P\Big(
\sum_{i\in J}\xi_i^{(n)}<\frac K2
\Big)\le e^{-cK}\le e^{-cN^\kappa},
\qquad J\in\{J^+,J_n^-\}.
\]
Thus, the probability that condition~(ii) fails
at any of the $k_0$ steps is at most
$O(Ne^{-cN^\kappa})$. Moreover,
Proposition~\ref{prop:restart-count} implies that
the probability of having fewer than $k_0/2$
effective steps is bounded by
$O((C/K)^{k_0/4})=O(e^{-cN^{2/3}})$. Combining these estimates,
we obtain
\[
\P(\mathcal G^c)
\lesssim N^{-D}+Ne^{-cN^\kappa}
+e^{-cN^{2/3}}
\lesssim N^{-D}.
\]

We next establish the contraction on $\mathcal G$.
If
\[
\|\delta^{(N+k_0)}\|_\infty\le 2(NK)^{-1/6},
\]
there is nothing to prove. Otherwise, by monotonicity,
$\|\delta^{(n+1)}\|_\infty>2(NK)^{-1/6}$
for every $N\le n<N+k_0$.
Proposition~\ref{prop:restart-step} then yields
\[
\|\delta^{(n+1)}\|_\infty
\le
\left(1-cN^{-2/3}K^{-1/3-4/\beta}\right)
\|\delta^{(n)}\|_\infty
\]
at every effective step, while
$\|\delta^{(n+1)}\|_\infty\le\|\delta^{(n)}\|_\infty$
always holds at the remaining steps.
Since at least $k_0/2$ steps are effective on $\mathcal G$,
the elementary inequality $(1-x)^m\le e^{-mx}$
for $0\le x\le1$, together with
\eqref{eq:restart-rigidity}, \eqref{eq:restart-exponent},
and $k_0\ge N^{2/3+\varepsilon}$, yields
\[
\begin{aligned}
\|\delta^{(N+k_0)}\|_\infty
&\le
\max\left\{
2(NK)^{-1/6},
e^{-ck_0N^{-2/3}K^{-1/3-4/\beta}/2}
\|\delta^{(n)}\|_\infty
\right\}\\
&\le
\max\left\{
3N^{-1/6-\kappa/6},
N^{-1/6+\kappa/6}e^{-cN^{2\varepsilon/3}}
\right\}\le N^{-1/6-\omega},
\end{aligned}
\]
for all sufficiently large $N$, where $c>0$ may
change from line to line and $\omega=\kappa/12$.
By monotonicity, the same bound holds
simultaneously for every $k_0\le k\le N$.
Together with the estimate on $\P(\mathcal G^c)$,
this proves \eqref{eq:restart}.

Finally, the upward transition preserves
the Gaussian $\beta$ ensemble row laws, so
$\widetilde x^{(N+k)}\sim\mu_{N+k,\beta}$.
Moreover, $\widetilde x^{(N+k)}$ depends only
on the independent starting row
$\widetilde x^{(N)}$ and the future innovations,
which are jointly independent of $\mathcal F_N$.
Hence $\widetilde x^{(N+k)}$ is independent
of $\mathcal F_N$, completing the proof.
\end{proof}

\section{Proof of the correlation estimate}\label{sec:correlation}
We prove Proposition~\ref{prop:increments} in this section. The main idea is using the the upward transition to show that the one step increments $X_{n+1}-X_n$ has a leading centered random term of order $N^{-1/3}$. Then summing over $L$ levels gives fluctuations of order $\sqrt{L/N^{2/3}}$, and the Doob maximal inequality provides uniform control over all intermediate levels. The proof is very similar to the counterpart of \cite{BCEHK-LFL}. 

In this section, we fix
$0<\tau<1/3$, which will be chosen small in terms of $\varepsilon$ in the
proof of Proposition~\ref{prop:correlation}. We first state the following one step estimate.

\begin{pro}[One step lower bound]\label{prop:increments}
There is $C>0$ such that the
following holds. With probability at least
$1-O(N^{-D})$,
\begin{equation}\label{eq:incrementlower}
 X_{n+1}-X_n\ge(n+1)^{-1/3}\bigl(\xi_1^{(n)}-1\bigr)
 -Cn^{-2/3+\tau}\bigl(1+\xi_1^{(n)}\bigr)
 \qquad\text{for all } \lceil N/2\rceil\le n<N,
\end{equation}
where $\xi_1^{(n)}$ is from \eqref{def.g-gamma}.
\end{pro}

 With the
above proposition in hand, we first use it to prove
Proposition~\ref{prop:correlation}.

\begin{proof}[Proof of Proposition~\ref{prop:correlation},
assuming Proposition~\ref{prop:increments}]
Recall $M=N-L$ with $1\le L\le N/2$. For $M\le r\le N$, set
\[
 S_r:=\sum_{n=M}^{r-1}(n+1)^{-1/3}
 (\xi_1^{(n)}-1).
\]
Since innovations $\xi_1^{(n)}\sim\operatorname{Gamma}(\beta/2,2/\beta)$, the centered random variable $\xi_1^{(n)}-1$ satisfies
\[
 \log\E e^{t(\xi_1^{(n)}-1)}
 \le Ct^2,\qquad |t|\le c_0.
\]
Since $n\ge N/2$ and $\{\xi^{(n)}_1\}_n$ are independent, Doob's exponential maximal
inequality yields
\begin{equation}\label{eq:martmax}
 \P\left(\max_{M\le r\le N}|S_r|>u\right)
 \le 2\exp\left[
 -c\min\left\{\frac{u^2 N^{2/3}}L,uN^{1/3}\right\}
 \right].
\end{equation}
Taking $u=N^{-1/3+\varepsilon/2}\sqrt{L}$, we obtain
\[
 \P\left(
 \max_{M\le r\le N}|S_r|
 >N^{-1/3+\varepsilon/2}\sqrt{L}
 \right)\lesssim N^{-D}.
\]
Moreover, the exponential tail of the Gamma
variables and a union bound give
\[
 \P\left(
 \max_{M\le n<N}\xi_1^{(n)}
 >N^{\varepsilon/4}
 \right)\lesssim N^{-D}.
\]
Thus, choosing $\tau<\varepsilon/4$ and summing
\eqref{eq:incrementlower}, we obtain, with
probability at least $1-O(N^{-D})$, uniformly
for all $M\le a<b\le N$,
\[
\begin{aligned}
 X_a-X_b&\le S_a-S_b+
 C\sum_{n=a}^{b-1}n^{-2/3+\tau}
 (1+\xi_1^{(n)})\\
 &\le 2\max_{M\le r\le N}|S_r|+C(b-a)N^{-2/3+\tau+\epsilon/4}
 \lesssim N^{\varepsilon/2}(N^{-1/3}\sqrt{L}+N^{-2/3}L).
\end{aligned}
\]
This proves \eqref{eq:correlation}.

Finally, fix $0<\kappa<2/3$ and set
$L=\lfloor N^{2/3-\kappa}\rfloor$.
Apply \eqref{eq:correlation} on $[N,N+L]$
with $\varepsilon=\kappa/4$.
Since
$
 (N+L)^{-2/3}L\le N^{-\kappa},
$
the corresponding threshold is $O(N^{-\kappa/4})$
for sufficiently large $N$. This proves
\eqref{eq:shortinterval}.
\end{proof}

In the remainder of this section, we prove
Proposition~\ref{prop:increments}. 

\begin{proof}[Proof of Proposition~\ref{prop:increments}]
Recall the quantiles $q_i^{(n)}$ defined in
\eqref{def.quantile}. For each $r$, let
\begin{equation}\label{eq:corr-rigidity}
 \mathcal R_n:=
 \left\{
 |x_i^{(n)}-q_i^{(n)}|
 \le n^{-1/6+\tau}\widehat i_n^{-1/3}
 \text{ for all }1\le i\le n
 \right\},
\end{equation}
where $\widehat i_n=\min\{i,n+1-i\}$ and $\tau$ is small enough.
The rigidity estimate \cite[Theorem~2.4]{BEY}
gives $\P(\mathcal R_n^c)\lesssim n^{-D}$
for every $D>0$.

Since $X_n=n^{1/6}(x_1^{(n)}-2\sqrt n)$,
we have the exact identity
\[
\begin{aligned}
X_{n+1}-X_n
&=(n+1)^{1/6}(x_1^{(n+1)}-x_1^{(n)})+\bigl((n+1)^{1/6}-n^{1/6}\bigr)
(x_1^{(n)}-2\sqrt n)\\
&-2(n+1)^{1/6}(\sqrt{n+1}-\sqrt n).
\end{aligned}
\]
On $\mathcal R_n$, rigidity gives
\[
\bigl((n+1)^{1/6}-n^{1/6}\bigr)
(x_1^{(n)}-2\sqrt n)=O(n^{-1+\tau}).
\]
Moreover, a simple calculation directly yields
\[
2(n+1)^{1/6}(\sqrt{n+1}-\sqrt n)
=(n+1)^{-1/3}+O(n^{-4/3}).
\]
Consequently, on $\mathcal R_n$, we have
\begin{equation}\label{eq:corr-gap}
 X_{n+1}-X_n
 =(n+1)^{1/6}(x_1^{(n+1)}-x_1^{(n)})
 -(n+1)^{-1/3}+O(n^{-1+\tau}).
\end{equation}
Thus, it suffices to obtain a lower bound on
$x_1^{(n+1)}-x_1^{(n)}$.

Evaluating the upward secular equation
\eqref{eq:secularup} at $x_1^{(n+1)}$ gives
\begin{equation}\label{eq:corr-upsecular}
 \frac{\xi_1^{(n)}}{
 \sqrt{n+1}(x_1^{(n+1)}-x_1^{(n)})}
 =\frac{x_1^{(n+1)}-g^{(n)}}{\sqrt{n+1}}
 -B_n,\qquad \mbox{where}~B_n:=\frac1{\sqrt{n+1}}\sum_{i=2}^n
 \frac{\xi_i^{(n)}}{x_1^{(n+1)}-x_i^{(n)}}.
\end{equation}
The key is to show that $B_n$ is close to one
from below.
To this end, let
\begin{equation}\label{def.Ir}
 I_n:=\{i:n^{3\tau}\le i\le n-n^{3\tau}\}.
\end{equation}
On $\mathcal R_n\cap\mathcal R_{n+1}$,
rigidity implies that
\begin{equation}\label{def.dn}
 0<x_1^{(n+1)}-x_i^{(n)}
 \le2\sqrt{n}-q_i^{(n)}+
 4n^{-1/6+\tau}:=d^{(n)}_i,
 \qquad i\in I_n.
\end{equation}
Since all terms in $B_n$ are positive by
interlacing property of particles, we obtain
\[
 B_n\ge
 \sqrt{\frac 1{n+1}}
 \sum_{i\in I_n}\frac{\xi_i^{(n)}}{d_i^{(n)}}.
\]
We next use the following claim to bound the estimates about coefficients $d^{(n)}_i$, whose proof is given after the proof of this proposition.

\begin{claim*}
For all sufficiently large $n$,
\begin{equation}\label{eq:bulksums}
 1-\frac1{\sqrt n}\sum_{i\in I_n}
 \frac1{d_i^{(n)}}
 \lesssim n^{-1/3+\tau},\qquad
 \max_{i\in I_n}
 \frac1{d_i^{(n)}}
 \lesssim n^{1/6-2\tau},\qquad
 \sum_{i\in I_r}
 \frac1{(d_i^{(n)})^2}
 \lesssim n^{1/3-\tau}.
\end{equation}
\end{claim*}
\medskip

By the above claim and Bernstein's inequality for independent
gamma variables, we have
\[
\begin{aligned}
\P\left(
 \frac1{\sqrt{n+1}}\sum_{i\in I_n}
 \frac{\xi_i^{(n)}}{d_i^{(n)}}
 <1-Cn^{-1/3+\tau}
 \right)&\le
\P\left(
 \frac1{\sqrt{n+1}}\sum_{i\in I_n}
 \frac{\xi_i^{(n)}-1}{d_i^{(n)}}
 <-Cn^{-1/3+\tau}
 \right)\\
&\le
2\exp\left[
-c\min\left\{
\frac{n^{1/6+2\tau}}
{\displaystyle\sum_{i\in I_n}
\frac1{(d_i^{(n)})^2}},
\frac{n^{-1/6+\tau}}
{\displaystyle\max_{i\in I_n}
\frac1{d_i^{(n)}}}
\right\}
\right]\lesssim e^{-cn^{3\tau}}.
\end{aligned}
\]
Consequently, on $\mathcal R_n\cap\mathcal R_{n+1}$
and outside an event of probability at most
$O(e^{-cn^{3\tau}})$,
\begin{equation}\label{eq:upbulk}
 B_n\ge1-Cn^{-1/3+\tau}.
\end{equation}
Moreover, the Gaussian tail of $g^{(n)}$ gives
$
 |g^{(n)}|\le n^{1/6+\tau}
$
outside an event of probability at most $O(e^{-cn^{3\tau}})$.
Combining this with \eqref{eq:upbulk},
$x_1^{(n+1)}\le2\sqrt{n+1}
+(n+1)^{-1/6+\tau}$,
and \eqref{eq:corr-upsecular}, we obtain
\[
 \frac{\xi_1^{(n)}}{
 \sqrt{n+1}(x_1^{(n+1)}-x_1^{(n)})}
 \le1+Cn^{-1/3+\tau}.
\]
Hence
\[
 x_1^{(n+1)}-x_1^{(n)}
 \ge
 \frac{\xi_1^{(n)}}{\sqrt{n+1}}
 \bigl(1-Cn^{-1/3+\tau}\bigr).
\]
Substituting this into \eqref{eq:corr-gap}
yields
\[
 X_{n+1}-X_n
 \ge(n+1)^{-1/3}(\xi_1^{(n)}-1)
 -Cn^{-2/3+\tau}(1+\xi_1^{(n)}),
\]
as desired.
Finally, a union bound proves that
\eqref{eq:incrementlower} holds simultaneously
for all such $n$ with probability
$1-O(N^{-D})$.
\end{proof}

It remains to prove the claim used in the above proof.

\begin{proof}[Proof of the claim]
By definition \eqref{def.dn}, we have 
$
 d_i^{(n)}\ge 2\sqrt n-q_i^{(n)}
 \gtrsim n^{-1/6}\hat i^{2/3}.
$
Hence, since $i\ge n^{3\tau}$ on $I_n$,
\[
 \max_{i\in I_n}\frac1{d_i^{(n)}}
 \lesssim n^{1/6-2\tau},
 \qquad
 \sum_{i\in I_n}\frac1{(d_i^{(n)})^2}
 \lesssim n^{1/3}\sum_{i\ge n^{3\tau}}i^{-4/3}
 \lesssim n^{1/3-\tau}.
\]
This proves the second and third estimates.
For the first estimate, since each interval
$(\gamma_{i+1}^{(n)},\gamma_i^{(n)})$
has semicircle mass $1/n$ and $(2-x)^{-1}$
is increasing, a Riemann-sum comparison gives
\[
\begin{aligned}
 \frac1{\sqrt n}\sum_{i\in I_n}
 \frac1{2\sqrt n-q_i^{(n)}}
 =\frac1n\sum_{i\in I_n}
 \frac1{2-\gamma_i^{(n)}}\ge
 \int_{-2}^2
 \frac{\rho_{\mathrm{sc}}(x)}{2-x}\,\dd x
 -Cn^{-1/3+\tau}\ge 1-Cn^{-1/3+\tau},
\end{aligned}
\]
where the omitted upper-edge contribution is
$O(n^{-1/3+\tau})$.

Finally, using $d^{(n)}_i-2\sqrt{n}+q^{(n)}_i=4n^{-1/6+\tau}$, which is from \eqref{def.dn} and $2\sqrt n-q_i^{(n)}
\asymp n^{-1/6}i^{2/3}$, we can obtain
\[
\begin{aligned}
 \frac1{\sqrt n}\sum_{i\in I_n}
 \left(
 \frac1{2\sqrt n-q_i^{(n)}}-\frac1{d_i^{(n)}}
 \right)\le
 \frac{4n^{-1/6+\tau}}{\sqrt n}
 \sum_{i\in I_n}
 \frac1{(2\sqrt n-q_i^{(n)})^2}
 \lesssim n^{-1/3}.
\end{aligned}
\]
Combining these estimates proves the first estimate, which completes the proof.
\end{proof}

\section*{Acknowledgement}
 Z.Bao and R. Xu are partially supported by Hong Kong RGC Grant GRF17304225 and 17303826.
\section*{AI disclosure}
All the  ideas belong to the authors, except for those in Appendix \ref{app:righttail} which were given by ChatGPT 6 Astra. This appendix extends the right tail upper bound result in \cite{BBBK-Tails} to the regime $0<\beta<1$. Without this extension, our results hold for all $\beta\geq 1$.  Although the main part of the paper does not use any idea provided by any AI models, we indeed had some discussion with ChatGPT 6 Astra during the preparation of the work. The interaction is described below:  By early September, the proof strategy of the correlation part, and also the coupling idea in Proposition \ref{prop:restart} were already clear for us. We nevertheless did not know how to prove Proposition \ref{prop:restart} at that time. From a conversation between one of us and ChatGPT 6 Astra, we obtained a proof strategy using energy method. However, we decided not to use this method and to develop an alternative one ourselves. This eventually leads to the current proof. Nevertheless, for interested readers, we provide ChatGPT’s energy method as a reference \cite{energymethod}, which is an  AI generated raw material, except for adjustments to the equation numbering and notations to match those of this paper. Finally, ChatGPT has also been used for language editing. 

\appendix 

\section{The right tail upper bound for \texorpdfstring{$0<\beta<1$}
{0<beta<1}}\label{app:righttail}
\begin{center}
    (by ChatGPT 6 Astra)
\end{center}

The tail estimates in
\cite[Theorems~1.1 and~1.2(i)]{BBBK-Tails}
provide all the bounds needed here, except for
the right-tail upper bound when $0<\beta<1$.
We establish this remaining case by localizing
the tridiagonal matrix into short intervals and
comparing each interval with a GOE matrix.

\begin{pro}\label{prop:smallbeta}
Fix $0<\beta<1$ and $0<\varepsilon<1$.
There exist $T,n_0>0$ such that, for $n\ge n_0$
and $T\le t\le n^{1/12}$,
\begin{equation}\label{eq:smallbetatail}
 \P\left(x_1^{(n)}\ge2\sqrt n+tn^{-1/6}\right)
 \le\exp\left[-\frac{2\beta}{3}
 (1-\varepsilon)t^{3/2}\right].
\end{equation}
\end{pro}

\begin{proof}
Let $A=\sqrt{\beta}H_{n,\beta}$ be the tridiagonal matrix
\[
A=
\begin{pmatrix}
g_1 & a_1 & & & \\
a_1 & g_2 & a_2 & & \\
& a_2 & \ddots & \ddots & \\
& & \ddots & g_{n-1} & a_{n-1}\\
& & & a_{n-1} & g_n
\end{pmatrix},
\]
where $g_i\sim\mathcal N(0,2)$ and
$a_i\sim\chi_{\beta(n-i)}$, $1\le i<n$,
are all independent.
This is the Dumitriu--Edelman tridiagonal model
for the Gaussian $\beta$-ensemble~\cite{DE}.
In particular, the eigenvalues of
$H_{n,\beta}=\beta^{-1/2}A$ have joint distribution
$\mu_{n,\beta}$ defined in \eqref{eq:rowlaw},
and hence
\[
x_1^{(n)}\overset{d}{=}
\beta^{-1/2}\lambda_{\max}(A).
\]

We use the following GOE tail estimate from
\cite[Theorem~1.2(i)]{BBBK-Tails}: for every $0<\eta<1$, there exists $U_\eta>0$
such that, for all sufficiently large $m$
and $U_\eta\le u\le m^{1/4}$,
\begin{equation}\label{eq:GOEtail}
 \P\left(
 \lambda_{\max}(H_{m,1})
 \ge2\sqrt m+um^{-1/6}
 \right)
 \le\exp\left[-\frac23(1-\eta)u^{3/2}\right],
 \qquad U_\eta\le u\le m^{1/4}.
\end{equation}

\medskip
\noindent\textit{Localization.}
For an integer $r\ge1$, define
\[
 \psi_k(i)=
 \begin{cases}
 \cos\bigl(\frac{\pi}{2}(i/r-k)\bigr),
       & |i/r-k|\le1,\\
 0, & \text{otherwise}.
 \end{cases}
\]
These functions satisfy
\[
 \sum_k\psi_k(i)^2=1,
 \qquad
 \sum_k(\psi_k(i+1)-\psi_k(i))^2\lesssim r^{-2}.
\]
Consequently, for every $v\in\mathbb R^n$,
\[
 \langle v,Av\rangle
 =\sum_k\langle\psi_kv,A\psi_kv\rangle+\sum_{i=1}^{n-1}a_iv_iv_{i+1}
 \sum_k(\psi_k(i+1)-\psi_k(i))^2.
\]
The supports of $\psi_k$ lie in intervals $I_k$
of length at most $2r$. Since
$\sum_i|v_iv_{i+1}|\le\|v\|_2^2$,
the variational principle yields
\begin{equation}\label{eq:localization}
 \lambda_{\max}(A)
 \le\max_k\lambda_{\max}(A_{I_k})
 +\frac{C\max_i a_i}{r^2}.
\end{equation}
Moreover, the standard chi-distribution tail bound
implies
\begin{equation}\label{eq:maxoffdiag}
 \P(\max_i a_i>2\sqrt n)\le ne^{-cn}.
\end{equation}
Thus, outside this rare event, the localization error
is at most $O(\sqrt n/r^2)$.

\medskip
\noindent\textit{Comparison with GOE.}
Consider a window
$I=\{s+1,\ldots,s+\ell\}$, where $\ell\le2r$.
Set
\begin{equation}\label{eq:comparisonsize}
 m_s=\left\lceil\beta(n-s)+(1-\beta)\ell\right\rceil.
\end{equation}
Since $n-s\ge\ell$, we have $m_s\ge\ell$ and
\[
 m_s-j\ge\beta(n-s-j),
 \qquad 1\le j<\ell.
\]
By the monotonicity of chi-distributions in their
degrees of freedom, we can couple $A_I$ with the
leading $\ell\times\ell$ principal block of $H_{m_s,1}$
so that their diagonal entries agree and the
off-diagonal entries of the latter are larger.
Since the largest eigenvalue of a tridiagonal
matrix with nonnegative off-diagonal entries
is increasing in those entries, we obtain
\begin{equation}\label{eq:windowcomparison}
 \P(\lambda_{\max}(A_I)\ge q)
 \le\P(\lambda_{\max}(H_{m_s,1})\ge q).
\end{equation}

Choose a small $\eta>0$ and then
$\delta>0$ sufficiently small depending on
$\beta,\eta$. Set
\begin{equation}\label{eq:windowchoice}
 r=\lceil\delta t n^{1/3}\rceil,
 \qquad
 q=\sqrt\beta(2\sqrt n+tn^{-1/6})
 -\frac{C\sqrt n}{r^2}.
\end{equation}
By \eqref{eq:comparisonsize},
\[
 m_s\le\beta n-\beta s+2(1-\beta)r+1.
\]
Using the concavity of the square root, we obtain
\begin{equation}\label{eq:windowgap}
 q-2\sqrt{m_s}
 \ge \sqrt\beta\,n^{-1/6}
 \left(t+\frac{s}{n^{1/3}}
 -C_\beta\delta t
 -\frac{C_\beta}{\delta^2t^2}\right)\ge \sqrt\beta\,n^{-1/6}
 \left((1-\eta)t+\frac{s}{n^{1/3}}\right).
\end{equation}
where the last inequality follows by choosing
$\delta$ sufficiently small and then $T$ sufficiently
large.

For windows with $s\le c_0n^{7/12}$, where $c_0>0$
is sufficiently small depending on $\beta$, we have
$m_s/(\beta n)=1+o(1)$ uniformly. Hence the
GOE-scaled threshold
$u_s=(q-2\sqrt{m_s})m_s^{1/6}$ satisfies
\[
 \beta^{2/3}(1-2\eta)
 \left(t+\frac{s}{n^{1/3}}\right)
 \le u_s
 \le C_\beta\left(t+\frac{s}{n^{1/3}}\right)
 \le m_s^{1/4}.
\]
The upper bound follows from
$m_s\ge\beta(n-s)$ and the choice of $c_0$.
Increasing $T$ also ensures $u_s\ge U_\eta$.
Therefore, \eqref{eq:GOEtail} and
\eqref{eq:windowcomparison} give
\begin{equation}\label{eq:windowsharp}
 \P(\lambda_{\max}(A_I)\ge q)
 \le\exp\left[
 -\frac{2\beta}{3}(1-C\eta)
 \left(t+\frac{s}{n^{1/3}}\right)^{3/2}
 \right].
\end{equation}

The windows have left endpoints $s=jr$,
with multiplicity at most two. By convexity,
\[
 \left(t+\frac{jr}{n^{1/3}}\right)^{3/2}
 \ge t^{3/2}+c\delta j t^{3/2}.
\]
Thus, summing \eqref{eq:windowsharp} over all
windows with $s\le c_0n^{7/12}$ gives
\[
 \sum_{s\le c_0n^{7/12}}
 \P(\lambda_{\max}(A_I)\ge q)
 \le C_{\beta,\eta}
 \exp\left[-\frac{2\beta}{3}
 (1-C\eta)t^{3/2}\right].
\]

For the remaining windows, \eqref{eq:windowgap}
implies
$
 q-2\sqrt{m_s}\ge cn^{1/12}.
$
We use the elementary bound
\begin{equation}\label{eq:GOEcoarse}
 \P\left(\lambda_{\max}(H_{m,1})
 \ge2\sqrt m+v\right)
 \le Cm e^{-cv^2},\qquad v\ge1.
\end{equation}
Indeed, it follows from
$\lambda_{\max}(H_{m,1})
\le\max_i d_i+2\max_i b_i$,
where $d_i$ and $b_i$ are the diagonal and
off-diagonal entries, together with Gaussian
and chi-distribution tails.
Since there are at most $n+2$ windows and
$m_s\le2n$, their total contribution is bounded by
$
 O(n^2e^{-cn^{1/6}}).
$

Combining the two window estimates with
\eqref{eq:localization} and \eqref{eq:maxoffdiag},
we obtain
\[
 \P\left(
 x_1^{(n)}\ge2\sqrt n+tn^{-1/6}
 \right)\lesssim
 \exp\left[-\frac{2\beta}{3}
 (1-C\eta)t^{3/2}\right]
 +n^2e^{-cn^{1/6}}+ne^{-cn}.
\]
Finally, choosing $\eta$ sufficiently small
relative to $\varepsilon$ and increasing $T,n_0$
absorbs the prefactor and the remaining errors,
since $t^{3/2}\le n^{1/8}$.
This proves \eqref{eq:smallbetatail}.
\end{proof}
\bibliographystyle{abbrv}
\bibliography{Ref}
\end{document}